\documentclass{amsart}
\usepackage[utf8]{inputenc}
\usepackage{a4wide}
\usepackage{euscript}
\usepackage{fullpage}
\usepackage[T2A]{fontenc}
\usepackage{amsfonts}
\usepackage{amssymb, amsthm}
\usepackage{amsmath}
\usepackage{mathtools}
\usepackage{graphicx}
\usepackage{geometry}
\usepackage{tikz}
\usepackage{bbm}

\numberwithin{equation}{section}

\theoremstyle{plain}
\newtheorem{thm}{Тheorem}
\newtheorem*{lm}{Lemma}
\newtheorem{lemma}{Lemma}
\newtheorem{st}{Statement}

\theoremstyle{definition}

\theoremstyle{remark}

\newcommand{\sg}{\sigma}

\numberwithin{thm}{section}
\numberwithin{lemma}{section}
\numberwithin{st}{section}
\numberwithin{cor}{section}

\title{Distribution of the Volume of Weighted Gaussian Simplex}
\author[T.~Moseeva]{Tatiana Moseeva}
\address{Tatiana Moseeva, Leonhard Euler International Mathematical Institute, Russia}
\email{polezina@yandex.ru}

\begin{document}

\thanks{The work was supported by the Foundation for the
Advancement of Theoretical Physics and Mathematics ``BASIS''}
\thanks{The work was supported by Ministry of Science and Higher Education of the Russian Federation, agreement № 075-15-2019-1619}

\begin{abstract}
Let $X_0, \ldots, X_l$ be independent standard Gaussian vectors in $\mathbb{R}^d$ such that $l \leqslant d$.  We derive an explicit formula for the distribution of the volume of weighted Gaussian simplex without the origin --- $l$-dimensional simplex $\mathrm{conv}(\sg_0X_0, \ldots, \sg_lX_l)$ ($\sg_0, \ldots, \sg_l > 0$). 
\end{abstract}

\maketitle

\section{Introduction and main result}
Let $X$ be a standard Gaussian vector in $\mathbb R^k$. By definition, its length is chi-distributed with $k$ degrees of freedom:
\begin{align*}
    |X| \stackrel{d}{=} \chi_k.
\end{align*}
This distributional identity has a  nice well-known generalization to several vectors. Namely, let $X_0,X_1, \ldots, X_l$ be random points sampled independently from the standard normal distribution on $\mathbb{R}^d$ ($l \leqslant d$). Then using some geometrical reasoning based on the <<height $\times$ base>> formula for the volume of a simplex and the symmetry property of the multidimensional Gaussian distribution (for details, see, e.g.,~\cite[Chapter~7]{tA03} or~\cite{NV14}) it can be shown that the volume of the $l$-dimensional  Gaussian simplex with the origin is distributed as follows:
\begin{align}\label{withorigin}
    |\mathrm{conv}(0, X_1, \ldots, X_l)| \stackrel{d}{=} \frac{1}{l!}\chi_{d-l+1}\dots\chi_d, 
\end{align}
where $|\cdot|$ denotes the $l$-dimensional volume. It follows immediately that for any $\sg_1, \ldots, \sg_l > 0$,
\begin{align*}
    |\mathrm{conv}(0, \sg_1X_1, \ldots, \sg_lX_l)| \stackrel{d}{=} \sg_1\dots\sg_l\frac{1}{l!}\chi_{d-l+1}\dots\chi_d.
\end{align*}

If we remove the origin from the set of the vertices, the situation becomes more subtle. Using the Blaschke--Petkantchin formula (see~Statement~\ref{BPf} below) Miles showed~\cite{Miles} that
\begin{align}\label{1908}
     \mathbb{E}|\mathrm{conv}(X_0, \ldots, X_l)|^p = \left[\frac{2^{l/2}\sqrt{l+1}}{l!}\right]^p\prod_{i = d - l + 1}^{d}{\frac{\Gamma((i+p)/2)}{\Gamma(i/2)}}.
\end{align}
From this, in view of
\begin{align*}
    \mathbb{E}\chi_k^p = 2^{p/2}\frac{\Gamma((k+p)/2)}{\Gamma(k/2)},
\end{align*}
he concluded (without a proof) that it follows from the method of moments that 
\begin{align*}
    |\mathrm{conv}(X_0, \ldots, X_l)| \stackrel{d}{=} \frac{\sqrt{l+1}}{l!}\chi_{d - l + 1}\dots\chi_d,
\end{align*}
see also~\cite[Theorem~2.5]{GKT19}.

Note that now the weighted case does not follows from the unweighted one. Specifically, in their study of the convex hulls of several Gaussian random walks, Randon-Furling and Zaporozhets~\cite[Theorem~6.1]{RFZ} derived as an auxiliary result the following formula:
\begin{align}\label{volmoments}
    \mathbb{E}|\mathrm{conv}(\sg_0X_0, \ldots, \sg_lX_l)|^p = \left[\frac{2^{l/2}\sg_0\dots\sg_l}{l!}\sqrt{\frac{1}{\sg_0^2} + \ldots + \frac{1}{\sg_l^2}}\right]^p\prod_{i = d - l + 1}^{d}{\frac{\Gamma((i+p)/2)}{\Gamma(i/2)}}.
\end{align}
Although~\eqref{1908} was already known, it took them three pages to prove it by using a generalization of the Sudakov--Tsirelson theorem obtained by Paouris and Pivivarov~\cite[Proposition~4.1]{PP13}. Again, it was suggested~\cite[Remark~6.2]{RFZ} that from the method of moments the distribution of the volume of the weighted Gaussian simplex could be obtained. Since for their main result the authors needed only the moments of the volume, they didn't provide a detailed proof. Also they suggested that it is possible to find a direct derivation avoiding the method of moments.

In this paper, we provide such a derivation, which is intuitive and geometric. Our main result is the following theorem.
\begin{thm}\label{res1}
Fix some $l = 1, \ldots, d$. Let  $X_0, \ldots, X_l$ be independent $d$-dimensional standard Gaussian vectors. Then for any $\sg_0, \ldots, \sg_l > 0$

\begin{align}\label{thm1}
    |\mathrm{conv}(\sg_0X_0, \ldots, \sg_lX_l)| \stackrel{d}{=} \frac{1}{l!}\sg_0\dots\sg_l\sqrt{\frac{1}{\sg_0^2} + \ldots + \frac{1}{\sg_l^2}}\chi_{d - l + 1}\dots\chi_d,
\end{align}
 where $\chi_{d-l+1}, \ldots, \chi_d$ are independent random variables such that for any $k = d - l + 1, \ldots, d$ the random variable $\chi_k$ has the chi distribution with $k$ degrees of freedom. 
\end{thm}



The paper is organized as follows. First we recall some necessary facts from stochastic geometry in Section \ref{facts}  and then present the proof of Theorem \ref{res1} in Section \ref{prf}.

\section{Auxiliary results}\label{facts}
To prove the main result  we will need the following statement (\cite{SW08}, Theorem 7.2.7):

\begin{st}[Blaschke-Petkantschin formula]\label{BPf}
Let $h:(\mathbb{R}^d)^{l+1} \rightarrow \mathbb R$ be non-negative measurable function, $l~\in \left\{1, . . . , d \right\}$. Then 
\begin{align*}
    &\int\limits_{(\mathbb{R}^d)^{l+1}}{f(x_0, \ldots , x_l) \mathrm{d}x_0 \ldots \mathrm{d}x_l} = \\ \notag
    &=(l!)^{d-l}b_{d,l}\int\limits_{A_{d,l}}{\int\limits_{E^{l+1}}{f(x_0, \ldots , x_l){|\mathrm{conv}(x_0, \ldots, x_l)|}^{d-l}\lambda_{E}(\mathrm{d}x_0)\ldots \lambda_E(\mathrm{d}x_l)}\mu_{d,l}(\mathrm{d}E)},
\end{align*}
where $b_{d,k} = \frac{\omega_{d-k+1} \cdots \omega_{d}}{\omega_1 \cdots \omega_k}$.
\end{st}

The following result of this section is probably well-known, however, for the reader's convenience we present it with proof.

\begin{lemma}\label{lm2}
Let  $\eta_1, \eta_2, \theta_1, \theta_2, \xi_1, \xi_2$ be positive with probability one random variables, such that

\begin{align} \label{eqlm2,1}
    \eta_1\cdot\theta_1 = \xi_1, \\ \label{eqlm2,2}
    \eta_2\cdot\theta_2 = \xi_2,
\end{align}
and $\eta_i$ и $\theta_i$ are independent for any $i \in \{1,2\}$. 
Let us assume that $\theta_1 \stackrel{d}{=} \theta_2$ and $\xi_1 \stackrel{d}{=} \xi_2$. Then $\eta_1 \stackrel{d}{=} \eta_2$.

\begin{proof}
Taking the logarithm of equations \eqref{eqlm2,1} and \eqref{eqlm2,2} we obtain
\begin{align*}
  \log{\eta_1} + \log{\theta_1} = \log{\xi_1}, \\
  \log{\eta_2} + \log{\theta_2} = \log{\xi_2}.
\end{align*}

Random variables in the left-hand side of both of these equations are independent thus the characteristic function of their sum is equal to the product of their characteristic functions. 
Therefore,

\begin{align*}
    \phi_{\log{\eta_1}}(t) \cdot \phi_{\log{\theta_1}}(t) = \phi_{\log{\xi_1}}(t)  = \phi_{\log{\xi_2}}(t) = \phi_{\log{\eta_2}}(t) \cdot \phi_{\log{\theta_2}}(t),
\end{align*}
where in the second step we used the fact that $\xi_1 \stackrel{d}{=} \xi_2$. Taking into account equation $\theta_1 \stackrel{d}{=} \theta_2$, we obtain $\phi_{\log{\eta_1}}(t) = \phi_{\log{\eta_1}}(t)$, and the lemma follows. 

\end{proof}

\end{lemma}

\section{Proof of Theorem \ref{res1}}\label{prf}

Let us consider the case $l = d$ first. 

By the simplex volume formula,
\begin{align*}
    |\mathrm{conv}(\sg_0X_0, \ldots, \sg_dX_d)| &= |\mathrm{conv}(0, \sg_1X_1 - \sg_0X_0, \ldots, \sg_dX_d - \sg_0X_0)|  \\ 
    &=
    \frac{1}{d!}\det [\sg_1X_1 - \sg_0X_0, \ldots, \sg_dX_d - \sg_0X_0]  \\ \\ &= \frac{1}{d!}\det 
    \begin{pmatrix}
    \sg_1X_1^{(1)} - \sg_0X_0^{(1)} &\ldots & \sg_dX_d^{(1)} - \sg_0X_0^{(1)} \\
    \vdots & \ddots & \vdots \\
    \sg_1X_1^{(d)} - \sg_0X_0^{(d)} & \ldots & \sg_dX_d^{(d)} - \sg_0X_0^{(d)}
    \end{pmatrix}  \\ \\ &=
    \frac{1}{d!}\det
    \begin{pmatrix}
    \sg_1X_1^{(1)} - \sg_0X_0^{(1)} &\ldots & \sg_1X_1^{(d)} - \sg_0X_0^{(d)} \\
    \vdots & \ddots & \vdots \\
    \sg_dX_d^{(1)} - \sg_0X_0^{(1)} & \ldots & \sg_dX_d^{(d)} - \sg_0X_0^{(d)}
    \end{pmatrix},
\end{align*}
where in the last step we used the fact that the determinant does not change after transposing. 

Note that columns of the obtained matrix are independent Gaussian vectors.  Let us denote the $k$-th column of this matrix by $Y_k$. 
$\mathbb{E}Y_k = 0$. 
\begin{align*}
    \mathrm{cov}(Y_k^{i}, Y_k^{j}) &= \mathrm{cov}(\sg_iX_i^{(k)} - \sg_0X_0^{(k)}, \sg_jX_j^{(k)} - \sg_0X_0^{(k)}) = 
    \sg_i\sg_j\mathrm{cov}(X_i^{(k)}, X_j^{(k)}) - \\
    &- \sg_i\sg_0\underbrace{\mathrm{cov}(X_i^{(k)}, X_0^{(k)})}_{=0} - \sg_j\sg_0\underbrace{\mathrm{cov}(X_0^{(k)}, X_j^{(k)})}_{=0} + \sg_0^2\mathrm{cov}(X_0^{(k)}, X_0^{(k)}) = \\ 
    &= \sg_i\sg_j\mathrm{cov}(X_i^{(k)}, X_j^{(k)}) + \sg_0^2 = \sg_i\sg_j\delta_{ij} + \sg_0^2, 
\end{align*}
where $\delta_{ij}$ is the Kronecker delta.

Therefore, 
$$\mathrm{cov}(Y_k) = 
\begin{pmatrix}
\sg_1^2+\sg_0^2 &\sg_0^2 &\ldots &\sg_0^2 \\
\sg_0^2 &\sg_2^2+\sg_0^2 &\ldots &\sg_0^2 \\
\vdots &\vdots &\ddots &\vdots \\
\sg_0^2 & \sg_0^2 &\ldots &\sg_d^2 + \sg_0^2
\end{pmatrix}.$$

Note that the covariance matrix of $Y_k$ does not depend on $k$. Let us denote this matrix by $\textbf{M}$. Then we can say, that $Y_k = \textbf{A}\overline{X_k}$, where $\overline{X_k}$ are independent standard Gaussian vectors and $\textbf{A}\textbf{A}^\mathrm{T} = \textbf{M}$. 

It follows that 
\begin{align}\label{l=d.mid}
    |\mathrm{conv}(\sg_0X_0, \ldots, \sg_dX_d)| = \frac{1}{d!}\det[\textbf{A}\overline{X_1}, \ldots, \textbf{A}\overline{X_d}] = \frac{1}{d!}\det \textbf{A}\det[\overline{X_1}, \ldots, \overline{X_d}]. 
\end{align}

First let us calculate the determinant of $\textbf{A}$.
Since $\textbf{A}\textbf{A}^\mathrm{T} = \textbf{M}$, we obtain $\det\textbf{A}  = (\det \textbf{M})^{\frac{1}{2}}$. 

\begin{gather*}
    \det \textbf{M} = \det \left( 
    \begin{pmatrix}
    \sg_1^2 &\ldots &0 \\
    \vdots &\ddots &\vdots \\
    0 &\ldots &\sg_d^2
    \end{pmatrix} + \sg_0^2\cdot I\right) = \\ = \sg_1^2\ldots\sg_d^2 + \sum_{k = 1}^d{\frac{\sg_1^2\ldots\sg_d^2}{\sg_k^2}\cdot\sg_0^2} =
    \sg_0^2\sg_1^2\ldots\sg_d^2\sum_{k = 0}^d{\frac{1}{\sg_k^2}}.
\end{gather*}

Consequently, $\det \textbf{A} = \sqrt{\det \textbf{M}} = \sg_0 \ldots \sg_d\sqrt{\frac{1}{\sg_0^2} + \ldots + \frac{1}{\sg_d^2}}$.

It remains to find the distribution of $\det[\overline{X_1}, \ldots, \overline{X_d}]$, that is, the volume of the parallelotope spanned by vectors $\overline{X_1}, \ldots, \overline{X_d}$. Let us denote that parallelotope by $P$. 

$$|P| = d!\cdot |\mathrm{conv}(0, \overline{X_1}, \ldots, \overline{X_d})|.$$ 
By \eqref{withorigin}, $|\mathrm{conv}(0, \overline{X_1}, \ldots, \overline{X_d})| \stackrel{d}{=} \frac{1}{d!}\chi_1 \dots \chi_d$, thus 
\begin{gather}\label{l=d.fin}
    \det[\overline{X_1}, \ldots, \overline{X_d}] \stackrel{d}{=} \chi_1\dots\chi_d.
\end{gather}

Substituting \eqref{l=d.fin} in \eqref{l=d.mid} we obtain 
\begin{align}\label{l=d}
    |\mathrm{conv}(\sg_0X_0, \ldots, \sg_dX_d)| \stackrel{d}{=} \frac{1}{d!}\sg_0 \ldots \sg_d\sqrt{\frac{1}{\sg_0^2} + \ldots + \frac{1}{\sg_d^2}}\chi_1\dots\chi_d,
\end{align}
where random variables $\chi_1, \ldots, \chi_d$ are independent, $\chi_k$ has the chi distribution with $k$ degrees of freedom for any $k = 1, \ldots, d$.   

Now suppose that $l < d$. 
Consider the affine hull
$$V_l := \mathrm{aff}(\sigma_0X_0, \ldots, \sigma_lX_l). $$
With probability one, $V_l$ is an affine $l$-plane. 
Denote by $O_{V_l}$ the orthogonal projection of the origin onto $V_l$. 
Consider the linear subspace 
$$W_l := V_l - O_{V_l}.$$

\begin{lm}
$W_l$ is uniformly distributed over the $l$-dimensional linear Grassmanian with respect to Haar measure and independently of $|\mathrm{conv}(\sg_0X_0, \ldots, \sg_lX_l)|$. 
\begin{proof}
Let $f:G_{d,l} \rightarrow \mathbb{R}^1$, $g: \mathbb{R}^1 \rightarrow \mathbb{R}^1$ be non-negative bounded measurable functions.
 
\begin{align*}
    &\mathbb{E}\left[f(W_l)g\left(|\mathrm{conv}(\sg_0X_0, \ldots, \sg_lX_l)|\right)\right] \\ 
    &= \int\limits_{(\mathbb{R}^d)^{l+1}}{f(\mathrm{aff}(\sg_0x_0, \ldots, \sg_lx_l)-O_{V_l})g(|\mathrm{conv}(\sg_0x_0, \ldots, \sg_lx_l)|)((2\pi)^{-d/2})^{l+1}e^{-\frac{|x_0^2|}{2}}\dots e^{-\frac{|x_l|^2}{2}} \mathrm{d}x_0\dots\mathrm{d}x_l}  \\
    &=\frac{1}{\prod_{i = 0}^l{\sg_i}}((2\pi)^{-d/2})^{l+1}
    \int\limits_{(\mathbb{R}^d)^{l+1}}{f(\mathrm{aff}(y_0, \ldots, y_l)-O_{V_l})g(|\mathrm{conv}(y_0, \ldots, y_l)|)e^{-\frac{|y_0|^2}{2\sg_0^2}}\dots  e^{-\frac{|y_l|^2}{2\sg_l^2}}\mathrm{d}y_0\dots\mathrm{d}y_l}.
\end{align*}

Using Statement \ref{BPf} yields

\begin{align*}
    &\mathbb{E}\left[f(W_l)g\left(|\mathrm{conv}(\sg_0X_0, \ldots, \sg_lX_l)|\right)\right] = \underbrace{(l!)^{d-l}b_{d,l}\frac{1}{\prod_{i = 0}^l{\sg_i}}((2\pi)^{-d/2})^{l+1}}_{=C}\cdot \\ 
    &\int\limits_{A_{d,l}}\int\limits_{E^{l+1}}{f(\mathrm{aff}(y_0, \ldots, y_l)-O_{V_l})g(|\mathrm{conv}(y_0, \ldots, y_l)|)e^{-\left(\frac{|y_0|^2}{2\sg_0^2}+ \ldots + \frac{|y_l|^2}{2\sg_l^2}\right)}|\mathrm{conv}(y_0, \ldots, y_l)|^{d-l}\prod_{i = 0}^l\lambda_E(\mathrm{d}y_i)\mu_{d,l}(\mathrm{d}E)} \\
    &= C\int\limits_{G_{d,l}}f(L)\int\limits_{L^{\bot}}\int\limits_{(L+a)^{l+1}}{g(|\mathrm{conv}(y_0, \ldots, y_l)|)e^{-\left(\frac{|y_0|^2}{2\sg_0^2}+ \ldots + \frac{|y_l|^2}{2\sg_l^2}\right)}|\mathrm{conv}(y_0, \ldots, y_l)|^{d-l}\prod_{i = 0}^l\lambda_{L+a}(\mathrm{d}y_i)\,\mathrm{d}a\,\nu_{d,l}(\mathrm{d}L)}.
\end{align*}
Consider the linear $l$-plane 
$$U := \mathrm{lin}(e_1, \ldots, e_l). $$
Fix some $L \in G_{d,l}$. 
Let $\mathbf{Q}$ be the orthogonal matrix such that $L = \mathbf{Q}U$. 
Then we have
\begin{align*}
    &\int\limits_{L^{\bot}}\int\limits_{(L+a)^{l+1}}{g(|\mathrm{conv}(y_0, \ldots, y_l)|)e^{-\left(\frac{|y_0|^2}{2\sg_0^2}+ \ldots + \frac{|y_l|^2}{2\sg_l^2}\right)}|\mathrm{conv}(y_0, \ldots, y_l)|^{d-l}\prod_{i = 0}^l\lambda_{L+a}(\mathrm{d}y_i)\,\mathrm{d}a} \\ 
    &= \int\limits_{U^{\bot}}\int\limits_{(U+a)^{l+1}}{g(|\mathrm{conv}(\mathbf{Q}y_0, \ldots, \mathbf{Q}y_l)|)e^{-\left(\frac{|\mathbf{Q}y_0|^2}{2\sg_0^2}+ \ldots + \frac{|\mathbf{Q}y_l|^2}{2\sg_l^2}\right)}|\mathrm{conv}(\mathbf{Q}y_0, \ldots,\mathbf{Q} y_l)|^{d-l}\prod_{i = 0}^l\lambda_{U + a}(\mathrm{d}y_i)\,\mathrm{d}a} \\
    &= \int\limits_{U^{\bot}}\int\limits_{(U+a)^{l+1}}{g(|\mathrm{conv}(y_0, \ldots, y_l)|)e^{-\left(\frac{|y_0|^2}{2\sg_0^2}+ \ldots + \frac{|y_l|^2}{2\sg_l^2}\right)}|\mathrm{conv}(y_0, \ldots, y_l)|^{d-l}\prod_{i = 0}^l\lambda_{U + a}(\mathrm{d}y_i)\,\mathrm{d}a}.
\end{align*}
Note that the obtained expression does not depend on $L$, denote it by $\mathcal{J}(g)$. 

Therefore,
\begin{align*}
    &\mathbb{E}\left[f(W_l)g\left(|\mathrm{conv}(\sg_0X_0, \ldots, \sg_lX_l)|\right)\right] = C\mathcal{J}(g)\int\limits_{G_{d,l}}f(L)\nu_{d,l}(\mathrm{d}L).
\end{align*}

Substituting $f = 1$ and $g = 1$ yields 
$$ 1 = C\mathcal{J}(1). $$

Let $A$ be some measurable subset of $G_{d,l}$. Consider $f = \mathbbm{1}_A, g = 1$. 
We have
$$\mathbb{P}[W_l \in A] = C\mathcal{J}(1)\int_A{\nu_{d,l}(\mathrm{d}L)} = \int_A{\nu_{d,l}(\mathrm{d}L)}, $$
therefore $W$ is uniformly distributed over $G_{d,l}$. 

Taking $f = 1, g = \mathbbm{1}_B$ for some measurable $B \in \mathbb{R}$ we obtain
$$\mathbb{P}[|\mathrm{conv}(\sg_0X_0, \ldots, \sg_lX_l)| \in B] = C\mathcal{J}(\mathbbm{1_B}).$$

Substituting $f = \mathbbm{1}_A$ and  $g = \mathbbm{1}_B$ finishes the proof. 

\end{proof}
\end{lm}

Let $P_l$ be the orthogonal projection from $\mathbb{R}^d$ onto the first $l$ coordinates, denote by $P_l^W$ its restriction to $W_l$.   

We have
\begin{align}
    \notag|\mathrm{conv}(\sigma_0P_lX_0, \ldots, \sigma_lP_lX_l)| &= 
    |\mathrm{conv}(\sigma_0P_lX_0 - P_lO_{V_l}, \ldots, \sigma_lP_lX_l - P_lO_{V_l})|  \\\notag &=
    |\mathrm{conv}(P_l^W(\sigma_0X_0 - O_{V_l}), \ldots, P_l^W(\sigma_0X_0 - O_{V_l})|  \\ &= 
    |\mathrm{conv}(\sigma_0X_0 - O_{V_l}, \ldots, \sigma_lX_l - O_{V_l})| \cdot |\mathrm{det}(P_l^W)| \\ \notag&=
    |\mathrm{conv}(\sigma_0X_0, \ldots, \sigma_lX_l)| \cdot |\mathrm{det}(P_l^W)|.
\end{align}

As noted above, the orthogonal projection of standard Gaussian vector onto the $l$-dimensional subspace has the standard Gaussian distribution in that subspace, therefore, using \eqref{l=d}, we have:

\begin{gather}\label{system1}
    |\mathrm{conv}(\sigma_0P_lX_0, \ldots, \sigma_lP_lX_l)| \stackrel{d}{=} \frac{1}{l!}\sigma_0\ldots\sigma_l
    \sqrt{\frac{1}{\sigma_0^2} + \dots + \frac{1}{\sg_l^2}}\chi_1\dots\chi_l,
\end{gather}
where $\chi_1, \ldots, \chi_l$ are independent chi-distributed random variables.  

Thus we have 
\begin{gather}\label{l1}
    |\mathrm{conv}(\sigma_0X_0, \ldots, \sigma_lX_l)| \cdot |\mathrm{det}(P_l^W)| \stackrel{d}{=} \frac{1}{l!}\sigma_0\ldots\sigma_l
    \sqrt{\frac{1}{\sigma_0^2} + \dots + \frac{1}{\sg_l^2}}\chi_1\dots\chi_l,
\end{gather} wherein $|\mathrm{conv}(\sigma_0X_0, \ldots, \sigma_lX_l)|$ and $|\mathrm{det}(P_l^W)|$ are independent.

Consider the linear subspace $W := \mathrm{lin}(0, Y_1, \ldots, Y_l)$, where $Y_i$ are independent standard Gaussian vectors in $\mathbb{R}^d$. 
With probability one, $W$ is a linear $l$-plane. Since $Y_i$ are independent and spherically invariant, $W$ is uniformly distributed over the $l$-dimensional Grassmanian. 

Thus,
\begin{gather}\label{last1}
    |\mathrm{conv}(0, Y_1, \ldots, Y_l)\cdot |\mathrm{det}(P_l^W)|
    =  |\mathrm{conv}(0, P_lY_1, \ldots, P_lY_l)|, 
\end{gather}
where $|\mathrm{det}(P_l^W)|$ does not depend on $|\mathrm{conv}(0, Y_1, \ldots, Y_l)|$. 

Since \eqref{withorigin},
\begin{gather}\label{last2}
    |\mathrm{conv}(0, Y_1, \ldots, Y_l)| \stackrel{d}{=} \frac{1}{l!}\chi_{d-l+1}\dots\chi_d,
\end{gather}
where $\chi_i$ are independent chi-distributed random variables.

Similarly, using the fact, that $P_lY_i$ are independent standard Gaussian vectors in $\mathbb{R}^l$, we obtain
\begin{gather}\label{last3}
    |\mathrm{conv}(0, P_lY_1, \ldots, P_lY_l)| \stackrel{d}{=}\frac{1}{l!}\chi_1\dots\chi_l. 
\end{gather}

Combining \eqref{last1}, \eqref{last2} and \eqref{last3}, we have 
\begin{gather}\label{l2}
    \frac{1}{l!}\chi_{d-l+1}\dots\chi_d \cdot |\mathrm{det}(P_l^W)|\stackrel{d}{=} \frac{1}{l!}\chi_1\dots\chi_l.
\end{gather}

Multiplying \eqref{l2} by $\sigma_0\ldots\sigma_l
    \sqrt{\frac{1}{\sigma_0^2} + \dots + \frac{1}{\sg_l^2}}$ we obtain
\begin{align}\label{l3}
    \frac{1}{l!}\sigma_0\ldots\sigma_l
    \sqrt{\frac{1}{\sigma_0^2} + \dots + \frac{1}{\sg_l^2}}\chi_{d-l+1}\dots\chi_d \cdot |\mathrm{det}(P_l^W)|\stackrel{d}{=}\frac{1}{l!}\sigma_0\ldots\sigma_l
    \sqrt{\frac{1}{\sigma_0^2} + \dots + \frac{1}{\sg_l^2}}\chi_1\dots\chi_l.
\end{align}

Applying Lemma~\ref{lm2} to \eqref{l1} and \eqref{l3} we conclude
 
\begin{gather*}
    |\mathrm{conv}(\sg_0X_0, \ldots, \sg_lX_l)| \stackrel{d}{=} \frac{1}{l!}\sg_0 \dots \sg_l\sqrt{\frac{1}{\sg_0^2} + \ldots + \frac{1}{\sg_l^2}} \;\chi_{d-l+1} \dots \chi_d,
\end{gather*}
which finishes the proof.

    \bibliographystyle{plain}
    \bibliography{literature}

\end{document}